\documentclass[letterpaper,10pt,conference]{IEEEtran}

\usepackage{graphicx}  % Written by David Carlisle and Sebastian Rahtz

\usepackage{bm, color}

\begin{document}

% paper title
\title{ A fast approach to Discontinuous Galerkin solvers for Boltzmann-Poisson transport systems 
for full electronic  bands and phonon scattering}
\author{\authorblockN{%
Irene M. Gamba\authorrefmark{1}, Armando Majorana\authorrefmark{2}, Jose A.
Morales\authorrefmark{1} and 
Chi-Wang Shu\authorrefmark{3}}
\authorblockA{\authorrefmark{1} ICES, The University of Texas at Austin, Austin, TX 78712} 
\authorblockA{\authorrefmark{2}Dipartimento di Matematica e Informatica,
Universit\`a di Catania, Catania, Italy}
\authorblockA{\authorrefmark{3}Division of Applied Mathematics, Brown
University, Providence, RI 02912}
e-mail: gamba@math.utexas.edu, majorana@dmi.unict.it,
jmorales@ices.utexas.edu, shu@dam.brown.edu
}
\maketitle
\begin{abstract}
\boldmath

The present work is motivated by the development of a fast DG based deterministic solver for
the extension  of the BTE  to a system of transport Boltzmann equations for full electronic
multi-band transport with intra-band scattering mechanisms.
Our proposed method allows to find scattering effects of high complexity, such as anisotropic electronic bands or full band computations, 
by simply using the standard routines of a suitable Monte Carlo approach only once. 
In this short paper, we restrict our presentation to the single band problem 
as it will be also valid in the multi-band system as well. 
We present preliminary numerical tests of this method using the 
Kane energy band model, for a 1-D 400nm $n^{+}-n-n^{+}$ silicon channel diode, showing moments
at $t=0.5$ps and $t=3.0$ps.

\end{abstract}

\section*{Introduction}

% GAMBA'S MODIF TO THE INTRO BELOW: USED FOR THE ABSTRACT ABOVE.
%The present manuscript is motivated by the development of a fast DG based deterministic solver for
%the extension  of the BTE  to a system of transport Boltzmann equations for full electronic
%multi-band transport with intra-band scattering mechanisms.
%Our proposed method mainly consists in finding scattering effects of high complexity (such as numerically
%computed full electronic bands) by simply using the standard routines of a suitable Monte
%Carlo approach only once. Our new proposed approach allows for an easy implementation of scattering terms having high complexity such as anisotropic electronic bands or full band calculations.
%In this short paper, we restrict our presentation to the single band problem as the
%proposed method described next will be also valid in the multi-band system as well. 

The semi-classical Boltzmann-Poisson system guarantees a good description of the dynamics of
electrons in modern semiconductor devices.
The equations of this model are given by
\begin{equation}\label{BE}
\frac{\partial f}{\partial t}
+ \frac{1}{\hbar} \nabla_{\mathbf{k}} \varepsilon \cdot \nabla_{\mathbf{x}} f -
 \frac{q}{\hbar} \mathbf{E} \cdot \nabla_{\mathbf{k}} f = Q(f).
\end{equation}
\begin{equation}\label{pois}
\nabla_{\mathbf{x}} \cdot \left[ \epsilon_{r}(\mathbf{x}) \mathbf{E}  \right] = \mbox{} -
 \frac{q}{\epsilon_{0}} \left[\rho(t,\mathbf{x})-N_D(\mathbf{x})\right] .
\end{equation}
In Eq.~(\ref{BE}), $f$ represents the electron probability density function (\emph{pdf}) in
phase space $\mathbf{k}$ at the physical location $\mathbf{x}$ and time $t$. $\mathbf{E}$ is
the electric field and $\varepsilon(\mathbf{k})$ is the energy-band function. 
The collision operator
$$
Q(f) = \int_{\Omega_{\mathbf{k}}} \! \! \left[ 
S(\mathbf{k}', \mathbf{k}) f(t, \mathbf{x}, \mathbf{k}') -
S(\mathbf{k}, \mathbf{k}') f(t, \mathbf{x}, \mathbf{k}) 
\right] d \mathbf{k}' 
$$
describes electron-phonon interactions through the kernel $S(\mathbf{k}', \mathbf{k})$. 
Physical constants $\hbar$ and $q$ are the Planck constant divided by $2\pi$ and the positive
electric charge, respectively. 
In Eq.~(\ref{pois}), $\epsilon_{0}$ is the dielectric constant in a vacuum,
$\epsilon_{r}(\mathbf{x})$ labels the relative dielectric function depending on the material,
$\rho(t,\mathbf{x})$ is the electron density, and $N_D(\mathbf{x})$ is the doping. 
The kinetic equation (\ref{BE}) is an equation in six dimensions (plus time if the device is
not in steady state) for a truly 3-D device.
This  high dimensionality has been a motivation for the BP system to be solved by the
Direct Simulation Monte Carlo (DSMC) methods~\cite{jaco89}. 
Yet we have proposed in~\cite{Cheng_08_CMAME} a deterministic approach based discontinuous
Galerkin (DG) method for solving Eqs.~(\ref{BE})-(\ref{pois}) that can be competitive. 
We refer to~\cite{Cheng_08_CMAME}   for a detailed description of DG and examples of applications
of the DG scheme to 1D diode and 2D double gate MOSFET devices. \\
\section*{The proposed method}
We assume that $\Omega_{\mathbf{k}}$ be a bounded domain of the $\mathbf{k}$-vector variable,
and we introduce a partition of it by means of a family of $N$ open cells                              
$C_{\alpha}$ such that, for every $\alpha$ and $\beta$,
$$
C_{\alpha} \subseteq \Omega_{\mathbf{k}} , \:
C_{\alpha} \cap C_{\beta} = \emptyset \: ( \alpha \neq \beta) , \:
\bigcup_{\alpha = 1}^{N} \overline{C_{\alpha}} = \Omega_{\mathbf{k}} \, .
$$
If we integrate the kinetic equation Eq.~(\ref{BE}) over the cell $C_{\alpha}$, 
then we obtain
\begin{eqnarray}
& & \hspace{-25pt}
\frac{\partial \mbox{ }}{\partial t} 
\int_{C_{\alpha}} \! \! f(t, \mathbf{x}, \mathbf{k}) \: d \mathbf{k}
+ \nabla_{\mathbf{x}} \cdot 
\int_{C_{\alpha}} \! \frac{1}{\hbar} \nabla_{\mathbf{k}} \varepsilon \,
f(t, \mathbf{x}, \mathbf{k}) \: d \mathbf{k}
\nonumber
\\
& & \hspace{-35pt} 
\mbox{} - \frac{q}{\hbar} \mathbf{E} \cdot
\int_{\partial C_{\alpha}} \! \! f(t, \mathbf{x}, \mathbf{k}) \, \mathbf{n} \: d \sigma
=
\int_{C_{\alpha}} \! \! Q(f) (t, \mathbf{x}, \mathbf{k}) \: d \mathbf{k} ,
\label{intkin}
\end{eqnarray}
where $\mathbf{n}$ is the normal to the surface $\partial C_{\alpha}$.

Any Galerkin method at the lowest order for the $\mathbf{k}$-vector variable, given by a piecewise 
constant approximation, assumes that \emph{in every cell $C_{\alpha}$ and for fixed
$\mathbf{x}$ and time $t$,  $f$ can be approximated by an unknown}
$f_{\alpha}(t, \mathbf{x})$.
This means that we are assuming $f$, for fixed $t$ and $\mathbf{x}$, to be constant in each
cell, except for the boundaries  of the cells, where $f$ is not even defined.
Physically, the unknown $f_{\alpha}(t, \mathbf{x})$ representing the approximated probability
density function of finding an electron at physical position $\mathbf{x}$ and time $t$, with
its wave-vector $\mathbf{k}$ belonging to the cell $C_{\alpha}$.
\\
Introducing the approximation for the distribution function $f$, we have
$$
\frac{\partial \mbox{ }}{\partial t} 
\int_{C_{\alpha}} \! \! f(t, \mathbf{x}, \mathbf{k}) \: d \mathbf{k}
\approx
M_{\alpha} \, \frac{\partial f_{\alpha}}{\partial t} (t, \mathbf{x}) ,
$$
where
$\displaystyle M_{\alpha} = \int_{C_{\alpha}} 1 \: d \mathbf{k}$
is the measure of the cell $C_{\alpha}$.
Now, if we define
\begin{eqnarray}
&&
\bm{\eta}_{\alpha} = \int_{C_{\alpha}} \! \frac{1}{\hbar} \nabla_{\mathbf{k}} \varepsilon 
\: d \mathbf{k} \, , \label{etak}
\\
&&
K_{\alpha \beta} = \int_{C_{\alpha}} d \mathbf{k} 
\int_{C_{\beta}} d \mathbf{k}' \, S(\mathbf{k}', \mathbf{k}) \, ,
\label{Kk}
\end{eqnarray}
then we have
$$
\nabla_{\mathbf{x}} \cdot \int_{C_{\alpha}} \! \frac{1}{\hbar} \nabla_{\mathbf{k}} 
\varepsilon \, f(t, \mathbf{x}, \mathbf{k}) \: d \mathbf{k}
\approx
\bm{\eta}_{\alpha} \cdot \nabla_{\mathbf{x}}  f_{\alpha}(t, \mathbf{x}) \, ,
$$
and
$$
\int_{C_{\alpha}} \hspace{-7pt} Q(f) (t, \mathbf{x}, \mathbf{k}) \: d \mathbf{k}
\approx 
\sum_{\beta=1}^{N} \! \left[ K_{\alpha \beta} \, f_{\beta}(t, \mathbf{x})
- K_{\beta \alpha} f_{\alpha}(t, \mathbf{x}) \right] \! .
$$
Therefore, we obtain a set of equations (for $\alpha = 1, 2, ..., N$),
which gives an approximation of the Boltzmann equation (\ref{BE})
$$
M_{\alpha} \, \frac{\partial f_{\alpha}}{\partial t} (t, \mathbf{x}) +
\bm{\eta}_{\alpha} \cdot \nabla_{\mathbf{x}}  f_{\alpha}(t, \mathbf{x})
- \frac{q}{\hbar} \mathbf{E} \cdot
\int_{\partial C_{\alpha}} \! \! f(t, \mathbf{x}, \mathbf{k}) \, \mathbf{n} \: d \sigma
$$
\begin{equation}
\mbox{} = 
\sum_{\beta=1}^{N} \! \left[ K_{\alpha \beta} \, f_{\beta}(t, \mathbf{x})
- K_{\beta \alpha} f_{\alpha}(t, \mathbf{x}) \right] .
\label{bte-dg} 
\end{equation}
Eqs.~(\ref{bte-dg}) contain yet the ``old'' unknown $f$ in the surface integral.
Here, $f$ can be approximated using $f_{\alpha}$ and other ``new'' unknows $f_{\gamma}$, where
the indexes $\gamma$ correspond to the nearest cells to $C_{\alpha}$.
The specific form of this transport term related to the electric field requires the use of some standard 
definition of the numerical flux $\hat{f}( f_{\alpha},f_{\gamma}... )$ according to the DG method 
(read the Appendix at the end of this paper for more details). After this step, Eqs.~(\ref{bte-dg}) become a set of $N$ partial differential equations in the
new $N$ unknowns $f_{\alpha}$.
We remark that the \emph{constant} coefficients $M_{\alpha}$, $\bm{\eta}_{\alpha}$ and
$K_{\alpha \beta}$ do not depend on the unknown $f$, but only on the domain decomposition, the
energy-band function $\varepsilon$ and the kernel $S(\mathbf{k}', \mathbf{k})$ of the
collision operator. 

The main difficulty in applying DG method to Eq.~(\ref{BE}) is to calculate the numerical
parameters $K_{\alpha \beta}$ for not simple analytical or real numerical band, as one tries
to apply quadrature formulas to Eq.~(\ref{Kk}).
\\
Here, we propose a very easy scheme to find the value of the parameters $K_{\alpha \beta}$ by
simply using the standard routines of a DSMC (Monte Carlo) solver {\bf only once} to determine the scattering process.
\\
To this aim, we consider the Boltzmann equation Eq.~(\ref{BE}), with zero electric field, for
spatially homogeneous solutions, i.e.
\begin{equation}
\frac{\partial f}{\partial t} = Q(f) \, .
 \label{eqBQ}
\end{equation}
We denote by $\Gamma(\mathbf{k})$ the total scattering rate
$$
\Gamma(\mathbf{k}) = 
\int_{\Omega_{\mathbf{k}}} S(\mathbf{k}, \mathbf{k}') \: d \mathbf{k}' .
$$
This is known for analytical band structure (for instance, the texbooks give its explicit
formulas for different materials) and it is used in DSMC code also in the full band case. 
Now, Eq.~(\ref{eqBQ}) writes
\begin{equation}
\frac{\partial f}{\partial t} = 
\int_{\Omega_{\mathbf{k}}} S(\mathbf{k}', \mathbf{k}) \, f(t, \mathbf{k}') \: d \mathbf{k}' - 
\Gamma(\mathbf{k}) \, f(t, \mathbf{k}) \, .
\label{eqB_om}
\end{equation}
Let be $\beta \in [1,N]$. Therefore, we define the initial data
$$
f(0, \mathbf{k}) = \phi(\mathbf{k}) =
\left\{
\begin{array}{ll}
1 & \mbox{ if } \mathbf{k} \in C_{\beta}
\\
0 & \mbox{ otherwise}
\end{array}
\right.
$$
Choose a small time step $\Delta t$ and \emph{solve} Eq.~(\ref{eqB_om}) using a DSMC procedure
{\bf only} in the small interval $[ 0, \Delta t ]$. 
So, we will know, with a reasonably good accuracy, the solution $f_{MC}(\Delta t, \mathbf{k})$
at time $\Delta t$.
Consider again Eq.~(\ref{eqB_om}). Since, 
\begin{eqnarray*}
&&
f(\Delta t, \mathbf{k}) \approx \phi(\mathbf{k}) + \Delta t
\\[5pt]
&& \mbox{ }
\times \left[
\int_{\Omega_{\mathbf{k}}} S(\mathbf{k}', \mathbf{k}) \, \phi(\mathbf{k}') \: d \mathbf{k}' - 
\Gamma(\mathbf{k}) \, \phi(\mathbf{k}) \right] ,
\end{eqnarray*}
we have
\begin{eqnarray}
&&
f(\Delta t, \mathbf{k}) \approx 
\left[ 1 - \Delta t \, \Gamma(\mathbf{k}) \right] \phi(\mathbf{k})
\nonumber \\[5pt]
&& \mbox{} \hspace{48pt} + 
\Delta t \int_{_{C_{\beta}}} S(\mathbf{k}', \mathbf{k}) \: d \mathbf{k}' \, .
\label{fdt}
\end{eqnarray}
Assuming the equivalence of BTE e DSMC, we replace Eq.~(\ref{fdt}) with
\begin{eqnarray}
&&
f_{MC}(\Delta t, \mathbf{k}) \approx 
\left[ 1 - \Delta t \, \Gamma(\mathbf{k}) \right] \phi(\mathbf{k})
\nonumber \\[5pt]
&& \mbox{} \hspace{48pt} + 
\Delta t \int_{_{C_{\beta}}} S(\mathbf{k}', \mathbf{k}) \: d \mathbf{k}' \, .
\label{fMCdt}
\end{eqnarray}
Now, $\phi(\mathbf{k})$ is the given initial data and we have found 
$f_{MC}(\Delta t, \mathbf{k})$ by means of DSMC; hence, Eq.~(\ref{fMCdt}) gives the parameter
$K_{\alpha \beta}$ :
\begin{eqnarray*}
&& 
\int_{C_{\alpha}} d \mathbf{k} 
\int_{C_{\beta}} d \mathbf{k}' \, S(\mathbf{k}', \mathbf{k}) 
\\
&& \mbox{} \approx
\frac{1}{\Delta t}
\int_{C_{\alpha}} \! \! \left[ f_{MC}(\Delta t, \mathbf{k}) - 
\left[ 1 - \Delta t \, \Gamma(\mathbf{k}) \right] \phi(\mathbf{k})
\right] d \mathbf{k}  
\\
&& \mbox{} =
\frac{1}{\Delta t} \left[
\int_{C_{\alpha}} f_{MC}(\Delta t, \mathbf{k}) \: d \mathbf{k} - 
\delta_{\alpha \beta} \, M_{\alpha} \right]
\\
&& \mbox{} +
\delta_{\alpha \beta} \int_{C_{\alpha}} \Gamma(\mathbf{k}) \: d \mathbf{k} ,
\end{eqnarray*}
where $\delta_{\alpha \beta} = 1$ if $\alpha = \beta$, and $0$ otherwise.
\\

The following table shows the errors between the \emph{exact} values of $K_{\alpha \beta}$ and
the values obtained by means of a DSMC code, when the Kane model for the energy band is used.
$$
\begin{array}{|c|c|c|}
 \hline
particles & \mbox{maximum error} & \mbox{mean value error}
\\ \hline \hline
10^{6}  & 0.06291  &   0.0044408 \\ \hline
10^{7}  & 0.01403  &   0.0014447 \\ \hline
10^{8}  & 0.00906  &   0.0004225 \\ \hline
5 \cdot 10^{8}  & 0.00271  &   0.0002024 \\ \hline
10^{9}  & 0.00145  &   0.0001313 \\ \hline \hline
\end{array}
$$
%Errors of the numerical coefficients  $K_{\alpha \beta}$ using DSMC with respect to
%\emph{exact} values.
%
This method does not rely on scattering term symmetries and this proposed new approach extends to 
highly complex scattering mechanisms such as anisotropic electronic bands or full band calculations.
\section*{Preliminary Numerical results}
For the one dimensional silicon $n^+-n-n^+$ $400$nm channel diode,
where the doping is ${N}_{D} =  5 \times 10^{17}$ cm$^{-3}$ in the
n$^+$  and ${N}_{D} = 2 \times 10^{15}$ $cm^{-3}$ in the n region,
we use $1440$ cells in $\mathbf{k}$-space and $N_{x}$ intervals in $\mathbf{x}$-space.
The applied potential is $2 \, V$.
We consider the Kane model for the energy band, and we show some quantities at time $t = 0.5 ps$
(a transient state) and $t = 3.0 ps$.
\begin{figure}[h]
\centering
\includegraphics[width=0.99\columnwidth, height=4cm]{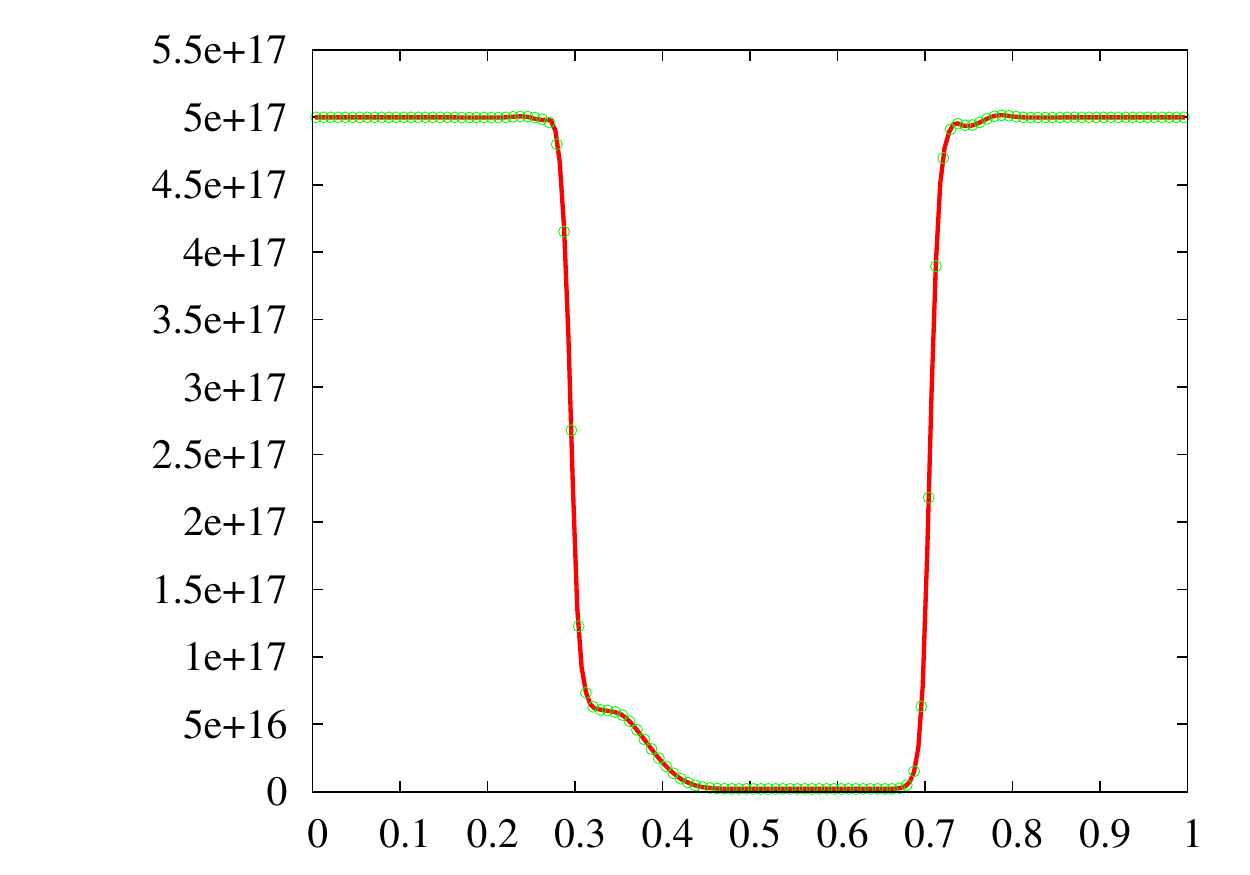}
\caption{Density of charge in $cm^{-3}$ at $t=0.5$ ps.
Continuous line ($N_{x} = 200$), points ($N_{x} = 120$).} 
\label{dens}
\end{figure}
\vspace{-10pt}
\begin{figure}[h]
\centering
\includegraphics[width=0.99\columnwidth, height=4cm]{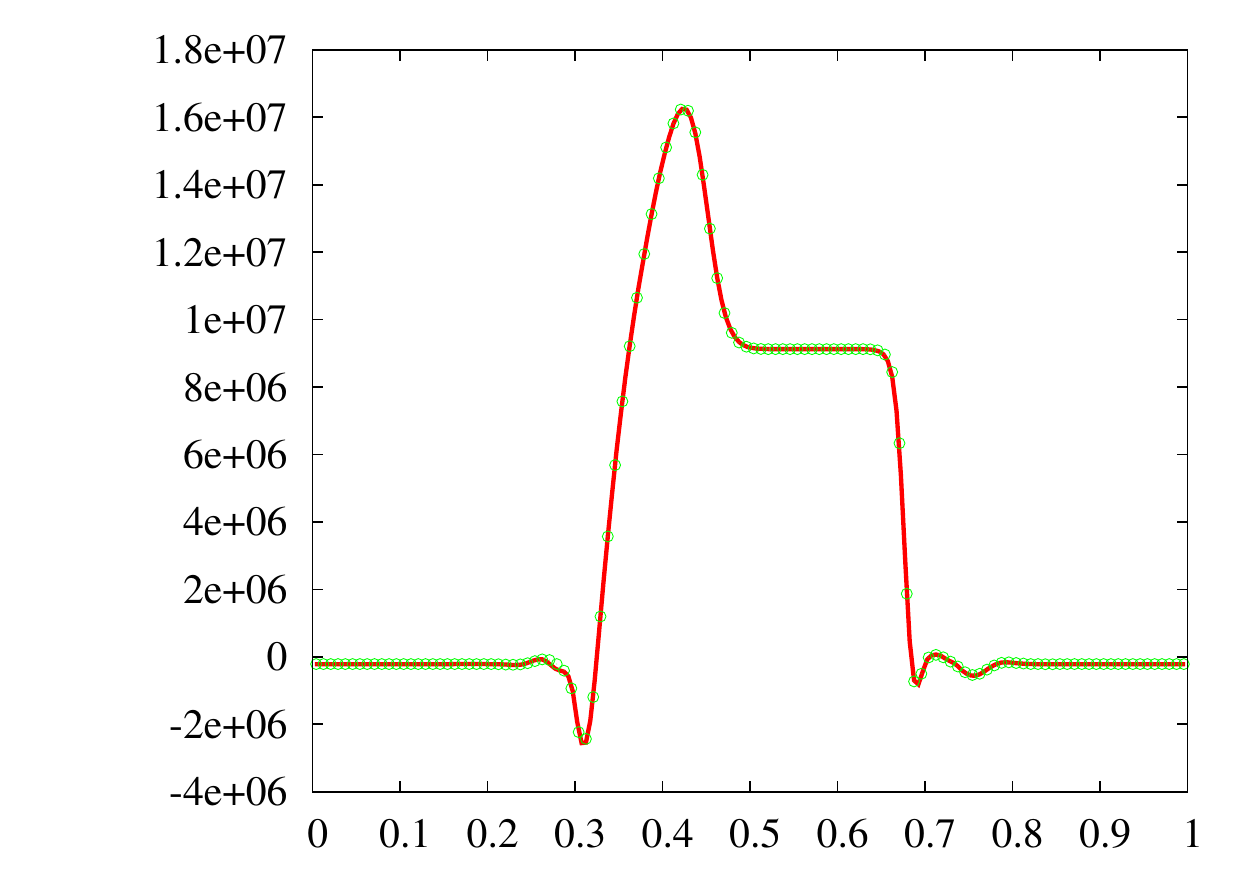}
\caption{Velocity in $cm/s$ at $t=0.5$ ps.
Continuous line ($N_{x} = 200$), points ($N_{x} = 120$).} 
\label{vel}
\end{figure}
\vspace{-10pt}
\begin{figure}[h]
\centering
\includegraphics[width=0.99\columnwidth, height=4cm]{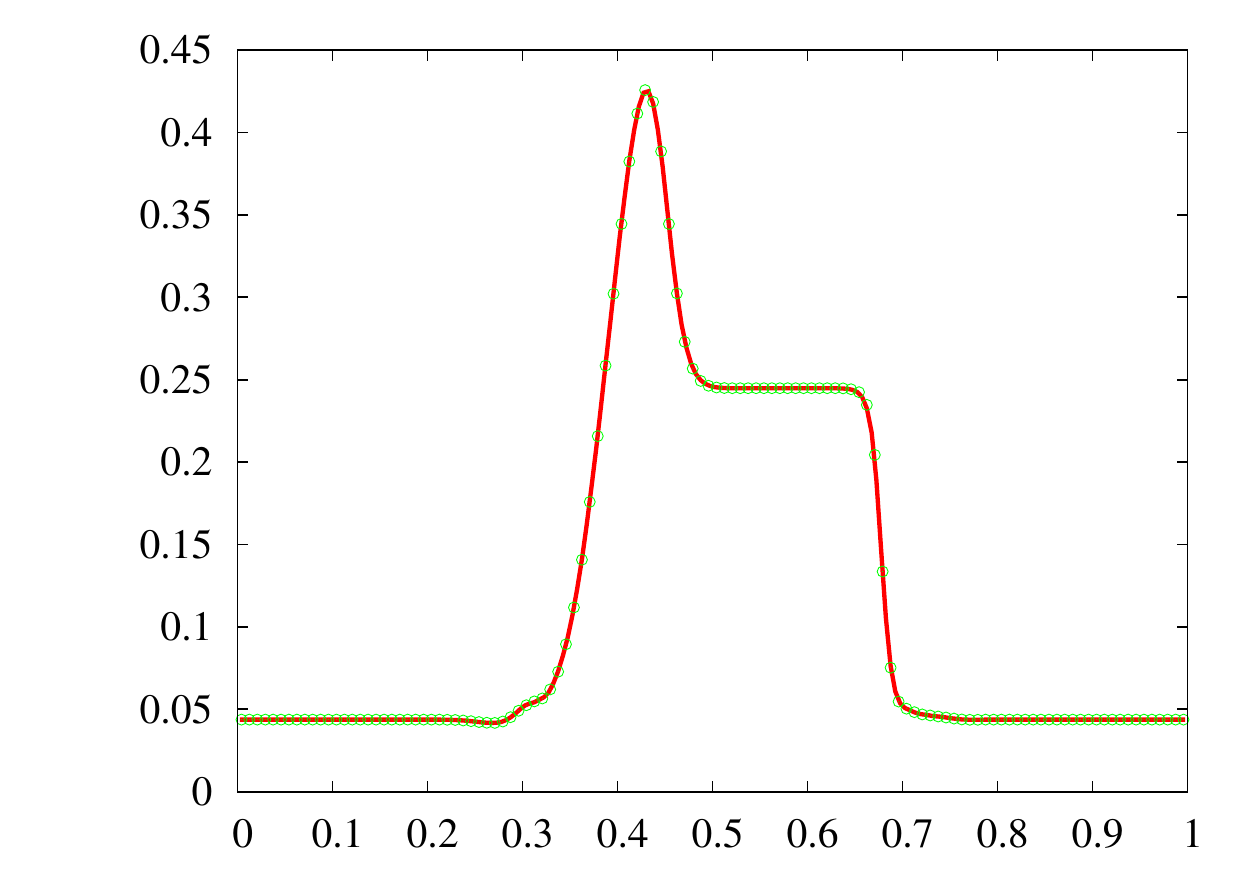}
\caption{Mean energy in $eV$ at $t=0.5$ ps.
Continuous line ($N_{x} = 200$), points ($N_{x} = 120$).} 
\label{ener}
\end{figure}
\begin{figure}[h]
\centering
\includegraphics[width=0.99\columnwidth, height=4cm]{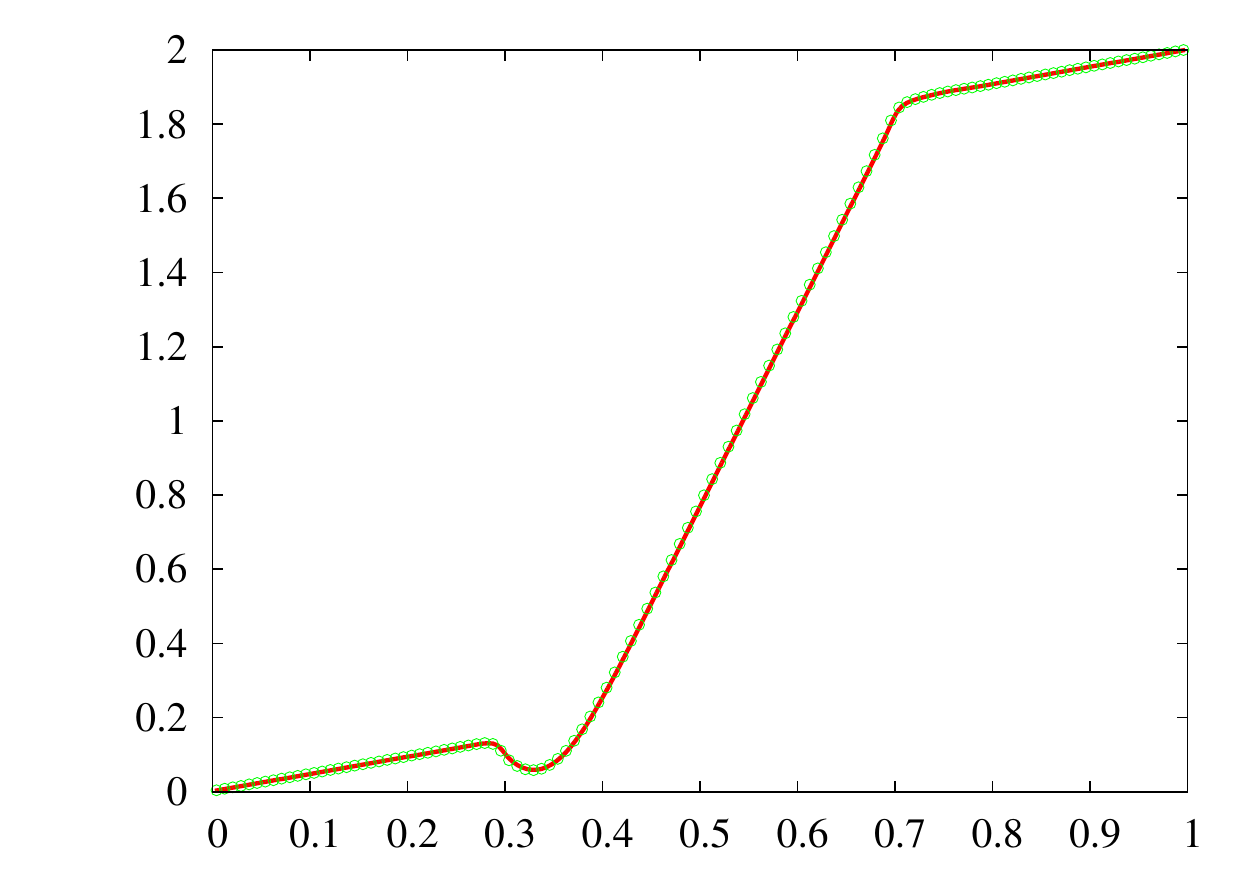}
\caption{Electric potential in $V$ at $t=0.5$ ps.
Continuous line ($N_{x} = 200$), points ($N_{x} = 120$).} 
\label{elec}
\end{figure}

$$
\begin{array}{|c|c|c|}
 \hline
nx & \mbox{minimum} & \mbox{maximum}
\\ \hline \hline
120 & 1.0517e-14  &  7.133e-04 \\ \hline
150 & 1.0517e-14  &  7.223e-04 \\ \hline
180 & 1.0520e-14  &  7.350e-04 \\ \hline
200 & 1.0519e-14  &  7.389e-04 \\ \hline \hline
\end{array}
$$
Minimum and maximum of $pdf$ multiplied a fixed function of $\mathbf{k}$ at $0.5$ps (a.u.)

\begin{figure}[h]
\centering
\includegraphics[width=0.99\columnwidth, height=4cm]{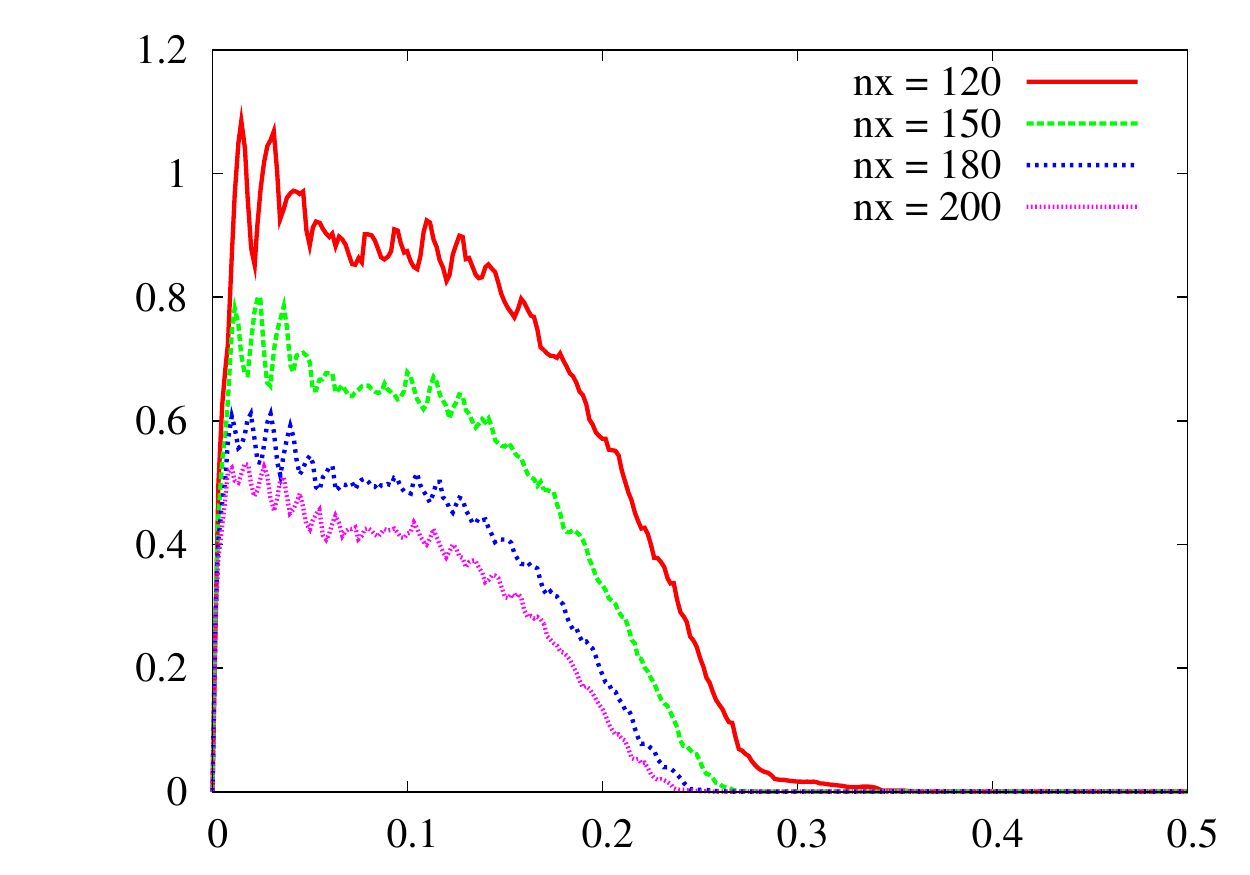}
\caption{The ratio, multiplies by $100$, of the number of cells in phase space where $pdf$ is
negative to the total numbers of cells versus time (in $ps$) for different $N_{x}$}
\end{figure}

%------------------------------------------

%
\begin{figure}[h]
\centering
\includegraphics[width=0.99\columnwidth, height=4cm]{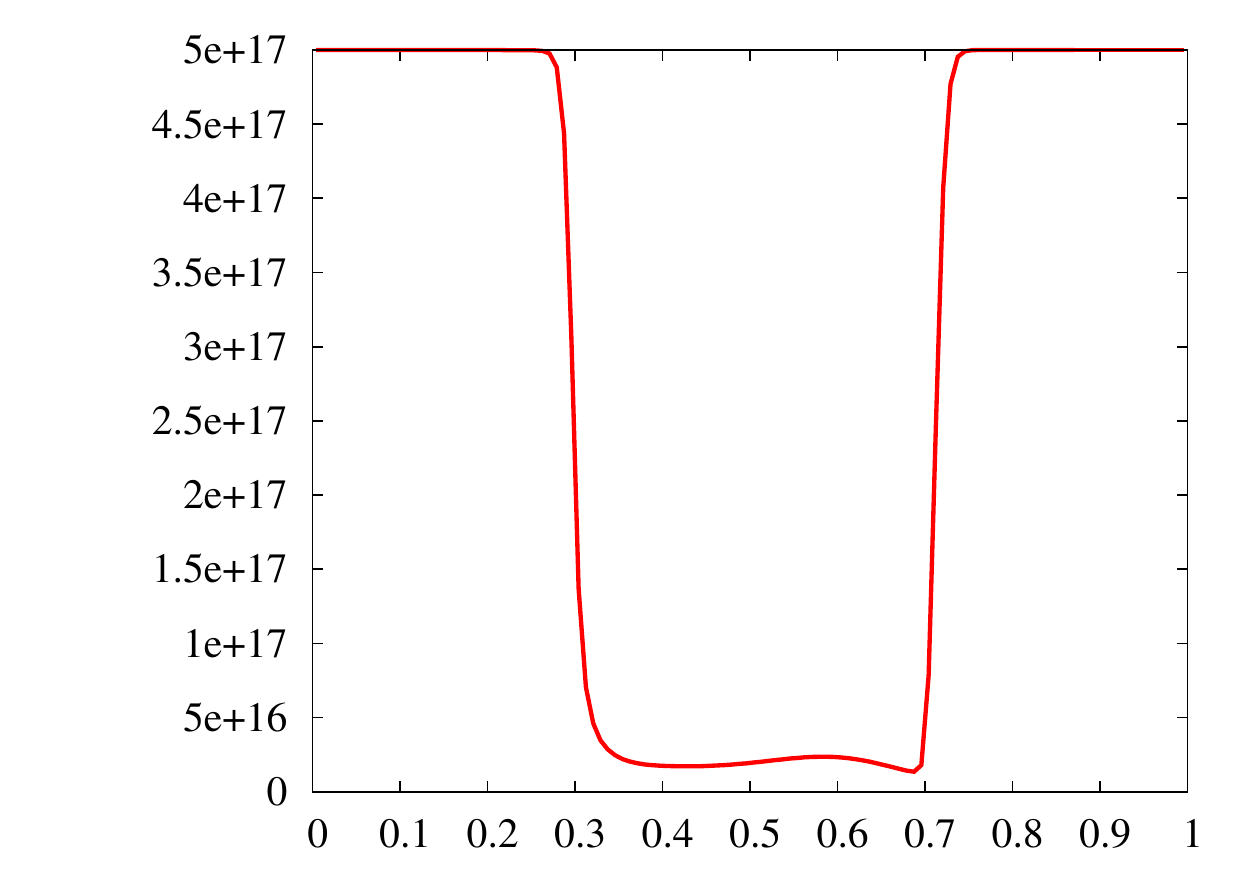}
\caption{Density of charge in $cm^{-3}$ at $t=3.0$ ps ($N_{x} = 120$).} 
\label{dens}
\end{figure}

\clearpage
\vspace{-10pt}
\begin{figure}[h]
\centering
\includegraphics[width=0.99\columnwidth, height=4cm]{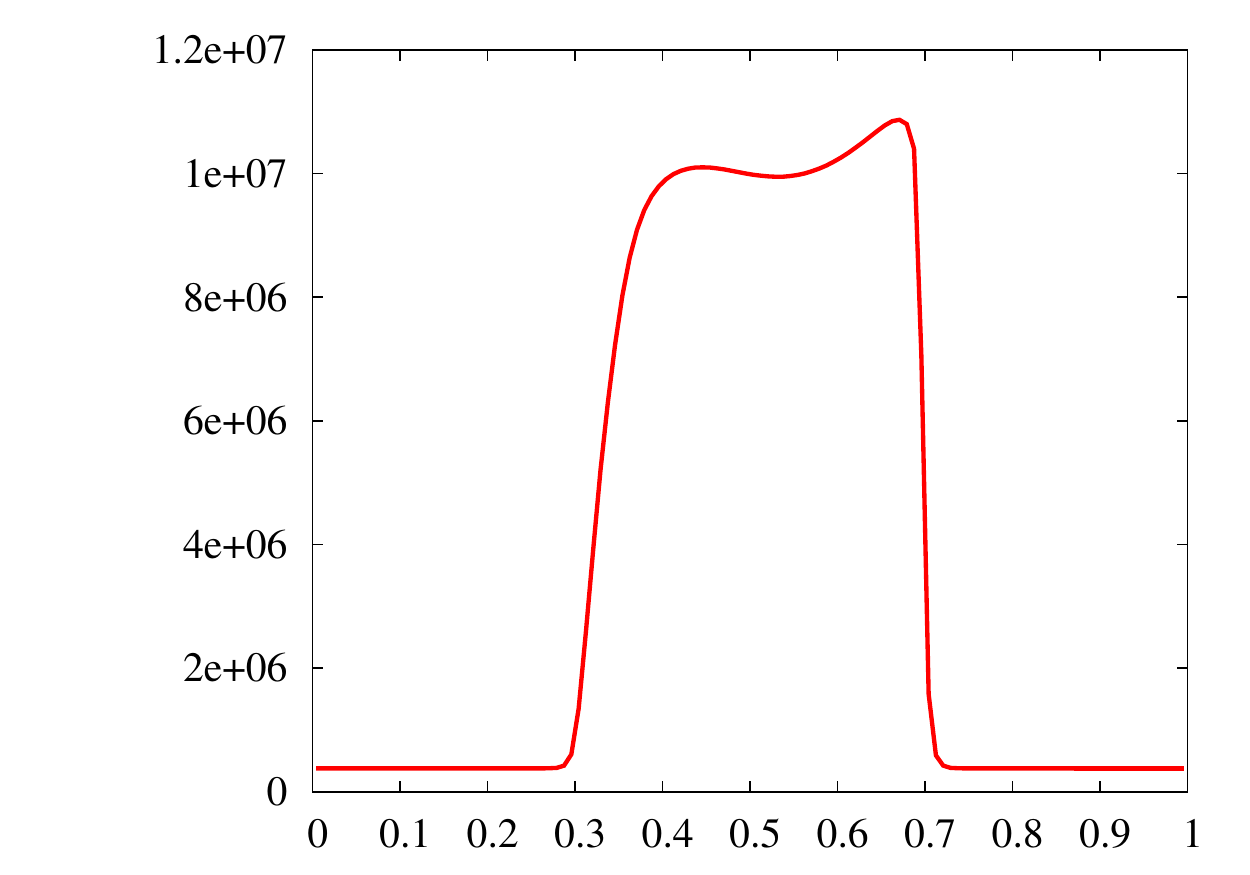}
\caption{Velocity in $cm/s$ at $t=3.0$ ps ($N_{x} = 120$).} 
\label{vel}
\end{figure}
\vspace{-10pt}
\begin{figure}[h]
\centering
\includegraphics[width=0.99\columnwidth, height=4cm]{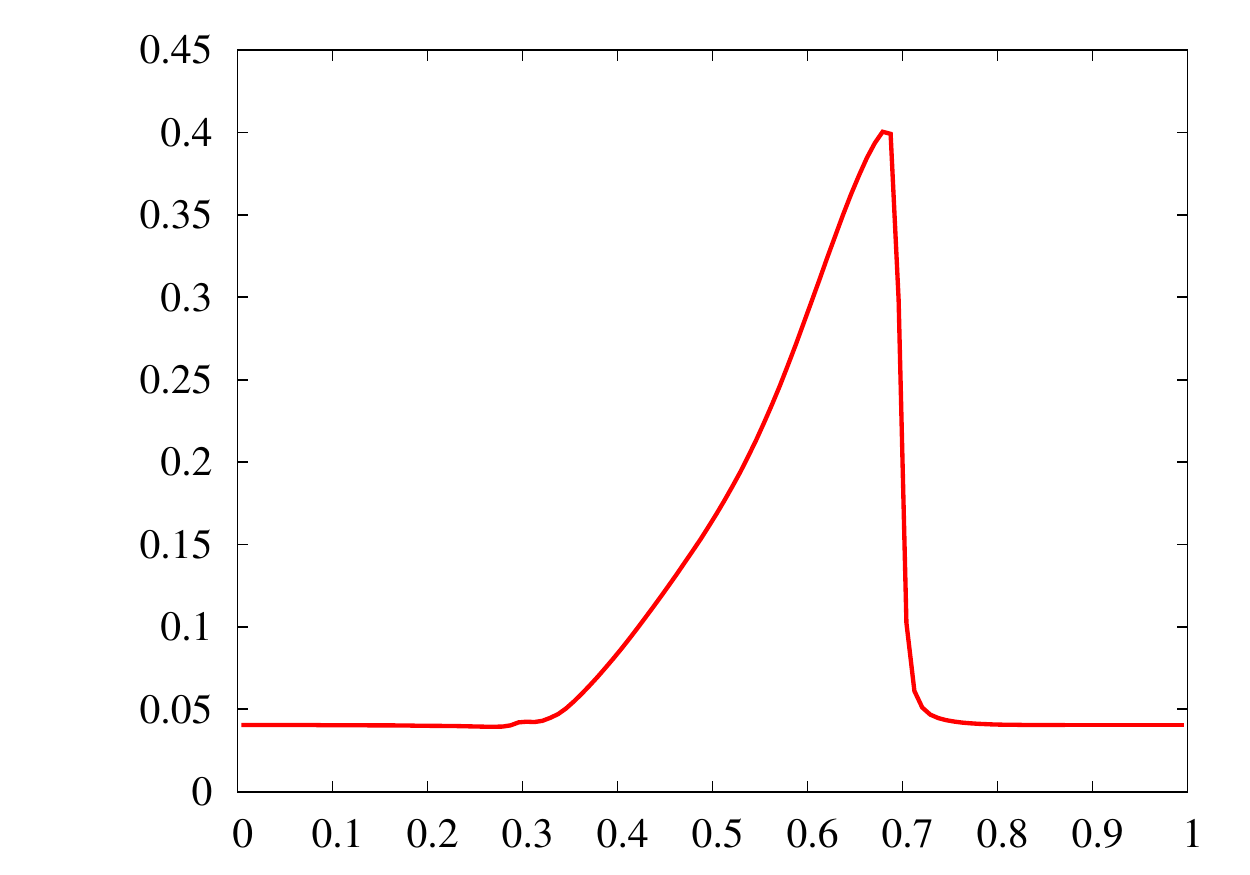}
\caption{Mean energy in $eV$ at $t=3.0$ ps ($N_{x} = 120$).} 
\label{ener}
\end{figure}
%

%----------------------------------------------

%\clearpage

% Appendix: Details about the Numerical Flux
% related to the transpor term Efield . grad_k f

\section*{Appendix: Treatment of the transport term related to the Electric Field}

%$ - \frac{q}{\hbar} \mathbf{E} \cdot \int_{\partial C_{\alpha}} \! \! f(t, \mathbf{x}, \mathbf{k}) \, \mathbf{n} \: d \sigma$}

The term $$ - \frac{q}{\hbar} \mathbf{E} \cdot 
\int_{\partial C_{\alpha}} \! \! f(t, \mathbf{x}, \mathbf{k}) \, \mathbf{n} \: d \sigma$$ still contains the original unknown pdf $f$, as it needs its value over the surface. However, this transport term related to the electric field $\mathbf{E}$ can be approximated by means of some standard definition of the Numerical Flux $\hat{f}$ according to the DG Method, adequate for a piecewise constant approximation. We will use the value $f_{\alpha}$ of the piecewise constant approximation for $f$ in the cell $C_{\alpha}$ and the values $f_{\gamma}$ in the nearest cells $C_{\gamma}$ neighboring $C_{\alpha}$. 

To illustrate this with a particular example, consider the case in which the electric field goes along the $x$ axis: $\mathbf{E} = (E_x,0,0)$. 
This 1D case has an associated cylindrical geometry in the $\mathbf{k}$-space:
$$
\mathbf{k} = k_* (u,r\cos\theta,r\sin\theta)
$$
where $k_*$ is a constant with dimensions of a $\mathbf{k}_{x_i}$-component, 
the normalized $u$ coordinate indicates the position along the $k_x$ axis, 
$r$ is the norm of the projection of the $\mathbf{k}$-point in a normalized $k_y$-$k_z$ plane, and $\theta\in[0,2\pi]$.

The particular symmetry of this case makes convenient to introduce annular $\mathbf{k}$-cells of the form 
$C_{\alpha} = [u_{a},u_{b}]\times[r_{a},r_{b}]\times[0,2\pi]$, related to the cylindrical geometry of the problem, and which look like rectangular cells on the $(u,r,\theta)$-space. 
Consider Figure \ref{urcell}, in which three neighboring $\mathbf{k}$-cells are shown: $C_{\alpha}$, $C_{\underline{\alpha}}$ (inferior to $C_{\alpha}$ in Fig. \ref{urcell}), and $C_{\overline{\alpha}}$ (superior to $C_{\alpha}$ in Fig. \ref{urcell}), as seen when projected in the $(u,r)$-space. 
Since in this case the 1D electric field is parallel to the $k_x$-axis, the transport term then reduces to:

\begin{figure}[h]
\centering
\includegraphics[width=0.70\columnwidth, height=5cm]{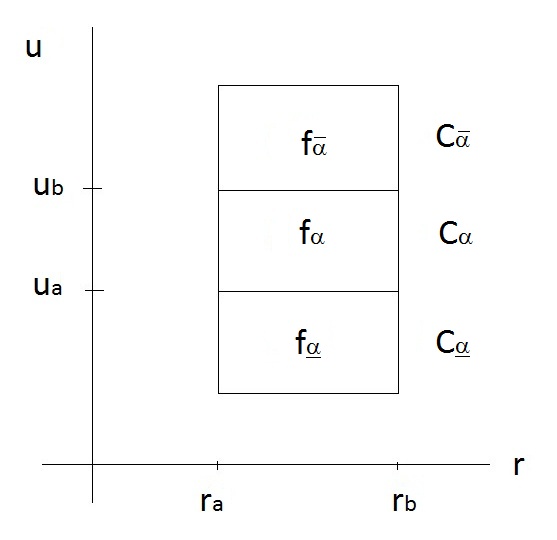}
\caption{Cell $C_{\alpha}$, and neighbor cells $C_{\underline{\alpha}}$ (inferior) and $C_{\overline{\alpha}}$ (superior), projected in the $(u,r)$-space}
\label{urcell}
\end{figure}

\begin{eqnarray*}
%&& 
 - \frac{q}{\hbar} \mathbf{E} \cdot \int_{\partial C_{\alpha}} f(t,\mathbf{x},\mathbf{k}) \, \mathbf{n} \: d \sigma \: = - \frac{q}{\hbar} E_x \int_{C_{\alpha}} \frac{\partial f}{\partial k_x} \: d \mathbf{k} &&\\
% &&
 \quad
  = - \frac{q}{\hbar} E_x \: k_{*}^{2} \int_{0}^{2\pi} d \theta \int_{r_a}^{r_b} d r \int_{u_a}^{u_b} d u \, \frac{\partial f}{\partial u} \: r  && \\
 %&&
 \quad
 = - \frac{q}{\hbar} E_x \: k_{*}^{2} \int_{0}^{2\pi} d \theta
\int_{r_a}^{r_b} r \: d r \left. \left[ \hat{f}(t,\mathbf{x},\mathbf{k}) \right] \right|_{u_{a}}^{u_{b}} & &
\end{eqnarray*}
For the transport term due to $-\mathbf{E}$ above, the Numerical Flux can be chosen according to the Upwind Principle. 
The flux over the considered boundaries is then:
\\
$$
\left. \hat{f}(t,\mathbf{x},\mathbf{k}) \right|_{u = u_{a}} = \left\{
\begin{array}{c l}
  f_{\alpha}(t,\mathbf{x}) & \mbox{if} \quad E_x \geq 0, \\
  f_{\underline{\alpha}}(t,\mathbf{x}) & \mbox{if} \quad E_x < 0.
\end{array}
\right.
$$

$$
\left. \hat{f}(t,\mathbf{x},\mathbf{k}) \right|_{u = u_{b}} = \left\{
\begin{array}{c l}
  f_{\overline{\alpha}}(t,\mathbf{x}) & \mbox{if} \quad E_x \geq 0, \\
  f_{\alpha}(t,\mathbf{x}) & \mbox{if} \quad E_x < 0.
\end{array}
\right.
$$
%
%
%optional acknowledgment
\section*{Acknowledgment}
The first and third authors are partially supported by
NSF DMS-1109525 and CHE-0934450. The second  is supported by PRIN 2009. 
The fourth author
is supported by NSF DMS 1112700 and DOE DE-FG02-
08ER25863. The second and fourth  author 
thanks the support from the J. Tinsley Oden Faculty
Research Fellowship from the Institute for Computational
Engineering and Sciences at the University of Texas at Austin.

%\clearpage
%\newpage

%put figures on the second page

%\begin{figure}
%\centering
%\includegraphics[width=\columnwidth]{fig1.eps}
%\caption{Sample Figure Caption [Times New Roman, 9pt]}
%\label{fig_sim}
%\end{figure}

% that's all folks
\end{document}